\documentclass[12pt]{amsart}

\usepackage{amssymb}
\usepackage{graphicx}
\usepackage{ifthen}
\usepackage{pict2e}
\usepackage{xargs}
\usepackage{xspace}

\newcommand{\definedterm}[1]{\emph{#1}}

\newcommand{\action}{\curvearrowright}
\newcommand{\Bairespace}[1][]{
  \ifthenelse{\equal{#1}{}}{\functions{\N}{\N}}{\functions{#1}{\N}}
}
\newcommand{\Bairetree}[1][]{
  \ifthenelse{\equal{#1}{}}{\functions{<\N}{\N}}{\functions{#1}{\N}}
}

\newcommand{\Borelchromaticnumber}[1]{\chi_B(#1)}

\newcommand{\calN}{\mathcal{N}}

\newcommandx{\CantorCantorspace}[2][1 =, 2 =]{
  \ifthenelse{\equal{#2}{}}
    {\ifthenelse{\equal{#1}{}}{\functions{\N \times \N}{2}}{\functions{#1 \times \N}{2}}}
    {\functions{#1}{(\functions{#2}{2})}}
}
\newcommand{\Cantorspace}[1][]{
  \ifthenelse{\equal{#1}{}}{\functions{\N}{2}}{\functions{#1}{2}}
}
\newcommand{\Cantortree}[1][]{
  \ifthenelse{\equal{#1}{}}{\functions{<\N}{2}}{\functions{<#1}{2}}
}
\newcommand{\cardinality}[1]{|#1|}
\newcommand{\closedopeninterval}[2]{[#1, #2)}
\newcommand{\closure}[1]{\overline{#1}}

\newcommand{\composition}{\circ}

\newcommandx{\concatenation}[2][1 = undefined, 2 = undefined]{
  \ifthenelse{\equal{#1}{undefined}}{{}\smallfrown}{
    \ifthenelse{\equal{#2}{undefined}}{\bigoplus #1}{\bigoplus_{#1} #2}
  }
}

\newcommandx{\constantfunction}[3][2 =, 3 =]{
  \ifthenelse{\equal{#2}{}}{c \from #1 \to \image{c}{#1}}{c_{#3} \from #1 \to #2}
}

\newcommandx{\convolution}[2][1 = undefined, 2 = undefined]{
  \ifthenelse{\equal{#1}{undefined}}{\mathrel{*}}{
    \ifthenelse{\equal{#2}{undefined}}{\bigotimes #1}{\bigotimes_{#1} #2}
  }
}
\newcommandx{\D}[2][2=]{\ifthenelse{\equal{#2}{}}{\mathbb{D}_{#1}}{\mathbb{D}_{#1,#2}}}
\newcommandx{\Deltaclass}[2][1=,2=]{
  \ifthenelse{\equal{#2}{}}{\mathbf{\Delta}_{#1}}{\mathbf{\Delta}^{#1}_{#2}}
}

\newcommand{\diameter}[2][]{\mathrm{diam}_{#1} \thinspace #2}
\newcommand{\differencefunction}[1]{\delta_{#1}}
\newcommand{\differenceset}[3][]{\Delta_{#1}(#2, #3)}

\newcommandx{\disjointunion}[2][1 =, 2 =]{
  \ifthenelse{\equal{#1}{}}{\sqcup}{
    \ifthenelse{\equal{#2}{}}{\bigsqcup #1}{{\bigsqcup_{#1} #2}}
  }
}

\newcommand{\domain}[1]{\mathrm{dom}(#1)}

\newcommand{\emptysequence}{\emptyset}
\newcommand{\equality}[1]{\Delta(#1)}
\newcommand{\equivalenceclass}[2]{[#1]_{#2}}
\newcommand{\Ethree}[1][]{\functions{\N}{\Ezero[#1]}}
\newcommand{\existnonmeagerlymany}{\exists^*}
\newcommand{\existsinfinitelymany}{\exists^\infty}

\newcommand{\extendedby}{\sqsubseteq}

\newcommand{\extensions}[2][]{
  \ifthenelse{\equal{#1}{}}{\calN_{#2}}{\calN_{#2} \intersection #1}
}
\newcommand{\Ezero}[1][]{\ifthenelse{\equal{#1}{}}{\mathbb{E}_0}{\mathbb{E}_{0, #1}}}

\newcommand{\flip}[1]{\overline{#1}}

\newcommand{\from}{\colon}

\newcommandx{\functions}[3][3 =]{
  \ifthenelse{\equal{#3}{}}{#2^{#1}}{#2^{#1}_{#3}}
}

\newcommand{\Gdelta}{$G_\delta$\xspace}

\newcommand{\goesto}{\rightarrow}

\newcommand{\Gzero}[1][]{\ifthenelse{\equal{#1}{}}{\mathbb{G}_0}{\mathbb{G}_{0, #1}}}

\newcommand{\GzeroN}[1][]{
  \ifthenelse{\equal{#1}{}}{\mathbb{G}_0^\N}{\mathbb{G}_{0, #1}^{\N}}
}

\newcommand{\identity}[1]{1_{#1}}
\newcommand{\image}[2]{#1(#2)}

\newcommandx{\identityfunction}[2][2 =]{
  \ifthenelse{\equal{#2}{}}{\mathrm{id} \from #1 \to #1}{\mathrm{id} \from #1 \to #2}
}

\newcommand{\infimum}[2][]{
  \ifthenelse{\equal{#1}{}}{\inf #1}{\inf_{#1}{#2}}
}
\newcommandx{\intersection}[2][1 =, 2 =]{
  \ifthenelse{\equal{#1}{}}{\cap}{
    \ifthenelse{\equal{#2}{}}{\bigcap #1}{{\bigcap_{#1} #2}}
  }
}
\newcommand{\inverse}[1]{#1^{-1}}
\newcommandx{\join}[2][1 =, 2 =]{
  \ifthenelse{\equal{#1}{}}{\vee}{
    \ifthenelse{\equal{#2}{}}{\bigvee #1}{{\bigvee_{#1} #2}}
  }
}
\newcommandx{\involutiongraph}[2][2=]{\ifthenelse{\equal{#2}{}}{\mathbb{G}_{#1}}{\mathbb{G}_{#1, #2}}}

\newcommand{\length}[1]{|#1|}
\newcommandx{\limit}[2][1 =, 2 =]{
  \ifthenelse{\equal{#1}{}}{\lim}{
    \ifthenelse{\equal{#2}{}}{\lim #1}{{\lim_{#1} #2}}
  }
}

\newcommand{\mathand}{\text{ and }}

\newcommand{\N}{\mathbb{N}}

\newcommand{\orbitequivalencerelation}[2]{E_{#1}^{#2}}
\newcommand{\orbitrelation}[2]{R_{#1}^{#2}}
\newcommand{\pair}[2]{(#1, #2)}

\newcommandx{\Piclass}[2][1=,2=]{
  \ifthenelse{\equal{#2}{}}{\mathbf{\Pi}_{#1}}{\mathbf{\Pi}^{#1}_{#2}}
}

\newcommand{\preimage}[2]{#1^{-1}(#2)}
\newcommandx{\product}[2][1 =, 2 =]{
  \ifthenelse{\equal{#1}{}}{\times}{
    \ifthenelse{\equal{#2}{}}{\prod #1}{{\prod_{#1} #2}}
  }
}
\newcommandx{\projection}[2][1 =, 2 =]{
  \ifthenelse{\equal{#1}{}}{\mathrm{proj}}{
    \ifthenelse{\equal{#2}{}}{\projection_{#1}}{
      \image{\projection[#1]}{#2}
    }
  }
}
\newcommand{\quadruple}[4]{(#1, #2, #3, #4)}

\renewcommandx{\restriction}[3][3 = undefined]{
  \ifthenelse{\equal{#3}{undefined}}{#1 \upharpoonright #2}{#1 \upharpoonright \pair{#2}{#3}}
}
\newcommandx{\sequence}[2][2 = undefined]{
  \ifthenelse{\equal{#2}{undefined}}{(#1)}{
    (#1)_{#2}
  }
}
\newcommandx{\set}[2][2 = undefined]{
  \ifthenelse{\equal{#2}{undefined}}{\{ #1 \}}{
    \{ #1 \suchthat #2 \}
  }
}
\newcommand{\setcomplement}[1]{\twiddle #1}
\newcommandx{\sets}[3][3 =]{
  \ifthenelse{\equal{#3}{}}{[#2]^{#1}}{[#2]^{#1}_{#3}}
}
\newcommandx{\Sigmaclass}[2][1=,2=]{
  \ifthenelse{\equal{#2}{}}{\mathbf{\Sigma}_{#1}}{\mathbf{\Sigma}^{#1}_{#2}}
}

\newcommand{\suchthat}{\mid}
\newcommand{\support}[1]{\text{supp}(#1)}

\newcommand{\twiddle}{\raisebox{1pt}{\scalebox{.75}{$\mathord{\sim}$}}}
\newcommandx{\union}[2][1 =, 2 =]{
  \ifthenelse{\equal{#1}{}}{\cup}{
    \ifthenelse{\equal{#2}{}}{\bigcup #1}{{\bigcup_{#1} #2}}
  }
}
\newcommand{\verticalsection}[2]{#1_{#2}}

\newcommand{\Baire}{Baire\xspace}
\newcommand{\Borel}{Bor\-el\xspace}
\newcommand{\Cauchy}{Cau\-chy\xspace}
\newcommand{\Hausdorff}{Haus\-dorff\xspace}
\newcommand{\Hjorth}{Hjorth\xspace}

\newcommand{\Kechris}{Kech\-ris\xspace}
\newcommand{\Klee}{Klee\xspace}
\newcommand{\Kuratowski}{Kur\-at\-ow\-ski\xspace}

\newcommand{\Lusin}{Lu\-sin\xspace}

\newcommand{\Montgomery}{Mont\-gom\-er\-y\xspace}
\newcommand{\Novikov}{No\-vik\-ov\xspace}
\newcommand{\Polish}{Po\-lish\xspace}

\newcommand{\Solecki}{Sol\-eck\-i\xspace}
\newcommand{\Todorcevic}{To\-dor\-cev\-ic\xspace}
\newcommand{\Ulam}{U\-lam\xspace}

\newenvironment{lemmaproof}{
  
  \begin{proof}
}{\end{proof}}

\newenvironment{propositionproof}{
  
  \begin{proof}
}{\end{proof}}

\newenvironment{theoremproof}{
  
  \begin{proof}
}{\end{proof}}

\newtheorem{lemma}{Lemma}[section]
\newtheorem{proposition}[lemma]{Proposition}

\newtheorem{theorem}[lemma]{Theorem}

\theoremstyle{definition}

\newenvironment{acknowledgements}{
  \textbf{Acknowledgements.}
}

\begin{document}


\begin{abstract}
  Under a mild definability assumption, we characterize the family of
  \Borel actions $\Gamma \action X$ of tsi \Polish groups on
  \Polish spaces that can be decomposed into countably-many
  actions admitting complete \Borel sets that are lacunary with
  respect to an open neighborhood of $\identity{\Gamma}$. In the
  special case that $\Gamma$ is non-archimedean, it follows that
  there is such a decomposition if and only if there is no continuous
  embedding of $\Ethree$ into $\orbitequivalencerelation{\Gamma}
  {X}$.
\end{abstract}

\author[B.D. Miller]{Benjamin D. Miller}

\address{
  Benjamin D. Miller \\
  Kurt G\"{o}del Research Center for Mathematical Logic \\
  Universit\"{a}t Wien \\
  W\"{a}hringer Stra{\ss}e 25 \\
  1090 Wien \\
  Austria
 }

\email{benjamin.miller@univie.ac.at}

\urladdr{
  http://www.logic.univie.ac.at/benjamin.miller
}

\thanks{The author was supported in part by FWF Grants P28153
  and P29999.}
  


\title[Essential countability]{Lacunary sets for actions of tsi groups}

\maketitle

\section*{Introduction}

The \definedterm{orbit equivalence relation} induced by a group
action $\Gamma \action X$ is the equivalence relation on $X$ given
by $x \mathrel{\orbitequivalencerelation{\Gamma}{X}} y \iff \exists
\gamma \in \Gamma \ \gamma \cdot x = y$. More generally, the
\definedterm{orbit relation} associated with a set $\Delta \subseteq
\Gamma$ is the binary relation on $X$ given by $x \mathrel
{\orbitrelation{\Delta}{X}} y \iff \exists \delta \in \Delta \ \delta \cdot x
= y$. A set $Y \subseteq X$ is \definedterm{$\Delta$-lacunary} if $y
\mathrel{\orbitrelation{\Delta}{X}} z \implies y = z$ for all $y, z \in Y$. 

Following the usual abuse of language, we say that an equivalence
relation $E$ on $X$ is \definedterm{countable} if $\cardinality
{\equivalenceclass{x}{E}} \le \aleph_0$ for all $x \in X$. We say that
a set $Y \subseteq X$ is \definedterm{$E$-complete} if
$\equivalenceclass{x}{E} \intersection Y \neq \emptyset$ for all $x
\in X$. The \definedterm{product} of equivalence relations $E_n$ on
$X_n$ is the equivalence relation $\product[n \in \N][E_n]$ on
$\product[n \in \N][X_n]$ given by $\sequence{x_n}[n \in \N] \mathrel
{(\product[n \in \N][E_n])} \sequence{y_n}[n \in \N] \iff \forall n \in \N
\ x_n \mathrel{E_n} y_n$. The \definedterm{$N$-fold power} of $E$ is
given by $\functions{N}{E} = \product[n \in N][E]$.

A \definedterm{graph} on $X$ is an irreflexive symmetric set $G
\subseteq X \times X$. We say that a set $Y \subseteq X$ is
\definedterm{$G$-independent} if $\restriction{G}{Y} = \emptyset$. A
\definedterm{$Z$-coloring} of $G$ is a map $\pi \from X \to Z$ such
that $\preimage{\pi}{\set{z}}$ is $G$-independent for all $z \in Z$.

A \definedterm{homomorphism} from a binary relation $R$ on $X$
to a binary relation $S$ on $Y$ is a map $\phi \from X \to Y$
such that $w \mathrel{R} x \implies \phi(w) \mathrel{S} \phi(x)$ for
all $w, x \in X$. More generally, a \definedterm{homomorphism} from
a sequence $\sequence{R_i}[i \in I]$ of binary relations on $X$ to a
sequence $\sequence{S_i}[i \in I]$ of binary relations on $Y$ is a
map $\phi \from X \to Y$ that is a homomorphism from $R_i$ to
$S_i$ for all $i \in I$. A \definedterm{reduction} of $R$ to $S$ is a
homomorphism from $\pair{R}{\setcomplement{R}}$ to $\pair{S}
{\setcomplement{S}}$, and an \definedterm{embedding} of $R$ into
$S$ is an injective reduction of $R$ to $S$.

Suppose that $\Gamma$ is a \Polish group and $X$ is a
\Borel space. We say that a \Borel action $\Gamma \action X$ is
\definedterm{$\sigma$-lacunary} if there are
$\orbitequivalencerelation{\Gamma}{X}$-invariant \Borel sets $X_n
\subseteq X$ with the property that $X = \union[n \in \N][X_n]$, open
neighborhoods $\Delta_n \subseteq \Gamma$ of $\identity{\Gamma}$,
and $\Delta_n$-lacunary $\orbitequivalencerelation{\Gamma}
{X_n}$-complete \Borel sets $B_n \subseteq X_n$ for all $n \in \N$.
A \Borel equivalence relation on a standard \Borel space is
\definedterm{essentially countable} if it is \Borel reducible to a
countable \Borel equivalence relation on a standard \Borel space.
The \Lusin-\Novikov uniformization theorem (see, for example, \cite
[Theorem 18.10]{Kechris}) easily implies that if $X$ is a standard \Borel
space, $\Gamma \action X$ is a $\sigma$-lacunary \Borel action, and
$\orbitequivalencerelation{\Gamma}{X}$ is \Borel, then
$\orbitequivalencerelation{\Gamma}{X}$ is essentially countable.

A well-known example of a non-essentially-countable \Borel
equivalence relation is the $\N$-fold power of the equivalence
relation $\Ezero$ on $\Cantorspace$ given by $c \mathrel{\Ezero}
d \iff \exists n \in \N \forall m \ge n \ c(m) = d(m)$.

A topological group is \definedterm{non-archimedean} if there is a
neighborhood basis of the identity consisting of open subgroups.
A topological group is \definedterm{tsi} if it has a compatible
two-sided-invariant metric. \Klee has shown that a \Hausdorff group
is tsi if and only if there is a neighborhood basis of the identity
consisting of conjugation-invariant open subsets (see \cite[1.5]
{Klee}). It follows that a \Hausdorff group is both non-archimedean
and tsi if and only if there is a neighborhood basis of the identity
consisting of normal open subgroups (see, for example, \cite[\S2]
{GaoXuan}).

\Hjorth-\Kechris have shown that if $\Gamma$ is a non-archimedean
tsi \Polish group, $X$ is a \Polish space, $\Gamma \action X$ is
\Borel, and $\orbitequivalencerelation{\Gamma}{X}$ is \Borel,
then either $\orbitequivalencerelation{\Gamma}{X}$ is essentially
countable or there is a continuous embedding of $\Ethree$ into
$\orbitequivalencerelation{\Gamma}{X}$ (see \cite[Theorem 8.1]
{HjorthKechris}). Our goal here is to give a classical proof of the
strengthening in which essential countability is replaced with
$\sigma$-lacunarity.

Given a graph $G$ on a \Borel space $X$, we write
$\Borelchromaticnumber{G} \le \aleph_0$ to indicate that $G$ has
\definedterm{countable \Borel chromatic number}, meaning that
there is a \Borel $\N$-coloring of $G$. \Kechris-\Solecki-\Todorcevic
have shown that there is a minimal analytic graph $\Gzero$ on a
standard \Borel space that does not have countable \Borel chromatic
number (see \cite[\S6]{KechrisSoleckiTodorcevic}).

In \S\ref{graph}, we characterize the class of increasing-in-$j$
sequences $\sequence{G_{i,j}}[i,j \in \N]$ of analytic graphs for
which there exist a function $f \from \N \to \N$ and a continuous
homomorphism $\phi \from \Cantorspace \to X$ from a sequence of
pairwise disjoint copies of $\Gzero$ to $\sequence{G_{i,f(i)}}[i \in
\N]$. In \S\ref{lacunary}, we show that for appropriately chosen
graphs, the inexistence of such homomorphisms yields
$\sigma$-lacunarity. In \S\ref{compositions}, we describe various
ways of refining such homomorphisms. And in \S\ref{dichotomies},
we establish a characterization of $\sigma$-lacunarity for \Borel actions
$\Gamma \action X$ of tsi \Polish groups with the property that
$\orbitrelation{\Delta}{X}$ is \Borel for every open set $\Delta \subseteq
\Gamma$. In the special case that $\Gamma$ is non-archimedean, this
yields our main result.

\section{A graph-theoretic dichotomy} \label{graph}

Fix $k_n \in \N$ such that $k_0 = 0$, $\forall n \in \N\ k_{n+1} \le \max
\set{k_m}[m \le n] + 1$, and $\forall k \in \N \existsinfinitelymany n
\in \N \ k_n = k$, as well as $s_n \in \Cantorspace[n]$ with the
property that $\forall k \in \N \forall s \in \Cantortree \exists n \in \N
\ (k = k_n \mathand s \extendedby s_n)$.

For all $s \in \Cantortree$, we use $\involutiongraph{s}$ to
denote the graph on $\Cantorspace$ given by $\involutiongraph
{s} = \set{\sequence{s \concatenation \sequence{i} \concatenation
c}[i < 2]}[c \in \Cantorspace]$. For all $k \in \N$, we use $\Gzero[k]$
to denote the graph on $\Cantorspace$ given by $\Gzero[k] = \union
[{\set{\involutiongraph{s_n}}[k = k_n \mathand n \in \N]}]$.

\begin{theorem} \label{graph:dichotomy}
  Suppose that $X$ is a \Hausdorff space and $\sequence{G_{i,j}}
  [i, j \in \N]$ is an increasing-in-$j$ sequence of analytic graphs on
  $X$. Then exactly one of the following holds:
  \begin{enumerate}
    \item There are \Borel sets $B_n \subseteq X$ such that $X =
      \union[n \in \N][B_n]$ and $\forall n \in \N \exists i \in \N \forall j
      \in \N \ \Borelchromaticnumber{\restriction{G_{i,j}}{B_n}} \le
      \aleph_0$.
    \item There exist a function $f \from \N \to \N$ and a continuous
      homomorphism $\phi \from \Cantorspace \to X$ from $\sequence
      {\Gzero[k]}[k \in \N]$ to $\sequence{G_{k, f(k)}}[k \in \N]$.
  \end{enumerate}
\end{theorem}

\begin{theoremproof}
  To see that conditions (1) and (2) are mutually exclusive, suppose
  that both hold, fix $n \in \N$ for which $\preimage{\phi}{B_n}$ is
  non-meager, fix $i \in \N$ such that $\forall j \in \N
  \ \Borelchromaticnumber{\restriction{G_{i,j}}{B_n}} \le \aleph_0$,
  fix a \Borel coloring $\psi \from B_n \to \N$ of $\restriction{G_{i,f(i)}}
  {B_n}$, fix $m \in \N$ for which $\image{(\inverse{\phi} \composition
  \inverse{\psi})}{\set{m}}$ is non-meager, fix $s \in \Cantortree$ for
  which $\image{(\inverse{\phi} \composition \inverse{\psi})}{\set{m}}$
  is comeager in $\extensions{s}$, and fix $\ell \in \N$ for which $i =
  k_\ell$ and $s \extendedby s_\ell$. It only remains to observe that
  there are comeagerly many $c \in \Cantorspace$ such that $s_\ell
  \concatenation \sequence{i} \concatenation c \in \image{(\inverse
  {\phi} \composition \inverse{\psi})}{\set{m}}$ for all $i < 2$,
  contradicting the fact that $\phi$ is a homomorphism from
  $\involutiongraph{s_\ell}$ to $G_{i, f(i)}$.
  
  It remains to show that at least one of conditions (1) and (2) holds.
  We can assume that $G_{i,j} \neq \emptyset$ for all $i,j \in \N$, in
  which case there are continuous surjections $\phi_{i, j} \from
  \Bairespace \to G_{i,j}$ for all $i, j \in \N$, as well as a continuous
  surjection $\phi_X \from \Bairespace \to \union[i,j \in \N][\image
  {\projection[X]}{G_{i,j}}]$.
  
  We will recursively define decreasing sequences $\sequence
  {X^\alpha_{i,j}}[\alpha < \omega_1]$ of subsets of $X$ such that
  $X^\alpha_{i,j} \subseteq X^\alpha_{i,j+1}$ and
  $\Borelchromaticnumber{\restriction{G_{i,j}}{\setcomplement
  {X^\alpha_{i,j}}}} \le \aleph_0$ for all $\alpha < \omega_1$ and $i, j
  \in \N$. We begin by setting $X^0_{i,j} = X$ for all $i, j \in \N$, and
  defining $X^\lambda_{i,j} = \intersection[\alpha < \lambda]
  [X^\alpha_{i,j}]$ for all $i, j \in \N$ and limit ordinals $\lambda <
  \omega_1$. To describe the construction of $X_{i,j}^{\alpha+1}$
  from $X_{i,j}^\alpha$, we require several preliminaries.
  
  We say that a quadruple $a = \quadruple{n^a}{f^a}{\phi^a}{\sequence
  {\psi^a_n}[{n < n^a}]}$ is an \definedterm{approximation} if $n^a \in
  \N$, $f^a \from \set{k_n}[n < n^a] \to \N$, $\phi^a \from
  \Cantorspace[n^a] \to \Bairespace[n^a]$, and $\psi^a_n \from
  \Cantorspace[n^a - 1 - n] \to \Bairespace[n^a]$ for all $n < n^a$.
  We say that an approximation $b$ is a \definedterm{one-step
  extension} of an approximation $a$ if:
  \begin{itemize}
    \item $n^a = n^b - 1$.
    \item $f^a = \restriction{f^b}{\set{k_n}[n < n^a]}$.
    \item $\forall i < 2 \forall s \in \Cantorspace[n^a] \ \phi^a(s)
      \extendedby \phi^b(s \concatenation \sequence{i})$.
    \item $\forall i < 2 \forall n < n^a \forall s \in \Cantorspace[n^a - n - 1]
      \ \psi_n^a(s) \extendedby \psi_n^b(s \concatenation \sequence{i})$.
  \end{itemize}
  We say that a quadruple $\gamma = \quadruple{n^\gamma}{f^\gamma}
  {\phi^\gamma}{\sequence{\psi^\gamma_n}[{n < n^\gamma}]}$ is a
  \definedterm{configuration} if $n^\gamma \in \N$, $f^\gamma \from
  \set{k_n}[n < n^\gamma] \to \N$, $\phi^\gamma \from \Cantorspace
  [n^\gamma] \to \Bairespace$, $\psi^\gamma_n \from \Cantorspace
  [n^\gamma - 1 - n] \to \Bairespace$ for all $n < n^\gamma$, and
  $(\phi_{k_n, f^\gamma(k_n)} \composition \psi^\gamma_n)(s) = \sequence
  {(\phi_X \composition \phi^\gamma)(s_n \concatenation \sequence
  {i} \concatenation s)}[i < 2]$ for all $n < n^\gamma$ and $s \in
  \Cantorspace[n^\gamma - n - 1]$. We say that a configuration
  $\gamma$ is \definedterm{compatible} with an approximation $a$
  if:
  \begin{itemize}
    \item $n^a = n^\gamma$.
    \item $f^a = f^\gamma$.
    \item $\forall s \in \Cantorspace[n^a] \ \phi^a(s) \extendedby
      \phi^\gamma(s)$.
    \item $\forall n < n^a \forall s \in \Cantorspace[n^a - n - 1]
      \ \psi^a_n(s) \extendedby \psi^\gamma_n(s)$.
  \end{itemize}
  We say that a configuration $\gamma$ is \definedterm{compatible}
  with a sequence $\sequence{X_{i,j}}[i, j \in \N]$ of subsets of $X$ if
  there is an extension $f \from \N \to \N$ of $f^\gamma$ with the
  property that $\image{(\phi_X \composition \phi^\gamma)}
  {\Cantorspace[n^\gamma]} \subseteq \intersection[i \in \N]
  [X_{i,f(i)}]$. We say that an approximation $a$ is \definedterm{$\sequence
  {X_{i,j}}[i,j \in \N]$-terminal} if no configuration is compatible with
  both a one-step extension of $a$ and $\sequence{X_{i,j}}[i,j \in \N]$.
  Let $A(a, \sequence{X_{i,j}}[i,j \in \N])$ denote the set of points of
  the form $(\phi_X \composition \phi^\gamma)(s_{n^a})$, where
  $\gamma$ varies over configurations compatible with both $a$
  and $\sequence{X_{i,j}}[i,j \in \N]$.
  
  \begin{lemma} \label{graph:dichotomy:splitting:one}
    Suppose that $\sequence{X_{i,j}}[i,j \in \N]$ is a sequence of
    subsets of $X$ and $a$ is an approximation for which $k_{n^a}
    \in \domain{f^a}$ and $A(a, \sequence{X_{i,j}}[i,j \in \N])$ is not
    $G_{k_{n^a},f^a(k_{n^a})}$-independent. Then $a$ is not
    $\sequence{X_{i,j}}[i,j \in \N]$-terminal.
  \end{lemma}
  
  \begin{lemmaproof}
    Fix configurations
    $\gamma_0$ and $\gamma_1$, compatible with $a$ and
    $\sequence{X_{i,j}}[i,j \in \N]$, for which $\sequence
    {(\phi_X \composition \phi^{\gamma_i})(s_{n^a})}[i < 2] \in
    G_{k_{n^a}, f^a(k_{n^a})}$. Then there exists $b \in \Bairespace$
    such that $\phi_{k_{n^a}, f^a(k_{n^a})}(b) = \sequence{(\phi_X
    \composition \phi^{\gamma_i})(s_{n^a})}[i < 2]$.
    Let $\gamma$ be the configuration given by $n^\gamma = n^a +
    1$, $f^\gamma = f^a$, $\phi^\gamma(s \concatenation \sequence
    {i}) = \phi^{\gamma_i}(s)$ for all $i < 2$ and $s \in \Cantorspace
    [n^a]$, $\psi^\gamma_n(s \concatenation \sequence{i}) = \psi^
    {\gamma_i}_n(s)$ for all $i < 2$, $n < n^a$, and $s \in
    \Cantorspace[n^a - n - 1]$, and $\psi^\gamma_{n^a}
    (\emptysequence) = b$. Then the unique approximation $b$ with
    which $\gamma$ is compatible is a one-step extension of $a$,
    so $a$ is not $\sequence{X_{i,j}}[i,j \in \N]$-terminal.
  \end{lemmaproof}
  
  \begin{lemma} \label{graph:dichotomy:splitting:two}
    Suppose that $\sequence{X_{i,j}}[i,j \in \N]$ is a sequence of
    subsets of $X$, $a$ is an approximation for which $k_{n^a}
    \notin \domain{f^a}$, and there exists $\ell \in \N$ such that
    $A(a, \sequence{X_{i,j}}[i,j \in \N])$ is not $G_{k_{n^a},
    \ell}$-independent. Then $a$ is not $\sequence{X_{i,j}}[i,j \in
    \N]$-terminal.
  \end{lemma}
  
  \begin{lemmaproof}
    Fix configurations $\gamma_0$ and $\gamma_1$, compatible
    with $a$ and $\sequence{X_{i,j}}[i,j \in \N]$, for which $\sequence
    {(\phi_X \composition \phi^{\gamma_i})(s_{n^a})}[i < 2] \in
    G_{k_{n^a}, \ell}$. By increasing $\ell$ if necessary, we can
    assume that $\image{\phi^{\gamma_0}}{\Cantorspace[n^a]}
    \union \image{\phi^{\gamma_1}}{\Cantorspace[n^a]} \subseteq
    X_{k_{n^a}, \ell}$. Fix $b \in \Bairespace$ such that
    $\phi_{k_{n^a}, \ell}(b) = \sequence{(\phi_X \composition
    \phi^{\gamma_i})(s_{n^a})}[i < 2]$, and let $\gamma$ be the
    configuration given by $n^\gamma = n^a + 1$, $f^\gamma(k) =
    f^a(k)$ for all $k < k_{n^a}$, $f^\gamma(k_{n^a}) = \ell$,
    $\phi^\gamma(s \concatenation \sequence{i}) = \phi^{\gamma_i}
    (s)$ for all $i < 2$ and $s \in \Cantorspace[n^a]$, $\psi^\gamma_n
    (s \concatenation \sequence{i}) = \psi^{\gamma_i}_n(s)$ for all $i
    < 2$, $n < n^a$, and $s \in \Cantorspace[n^a - n - 1]$, and
    $\psi^\gamma_{n^a}(\emptysequence) = b$. Then the unique
    approximation $b$ with which $\gamma$ is compatible is a
    one-step extension of $a$, so $a$ is not $\sequence{X_{i,j}}[i,j \in
    \N]$-terminal.
  \end{lemmaproof}
  
  As \Lusin's separation theorem (see, for example, \cite[Theorem
  14.7]{Kechris}) easily implies that every $G_{i,j}$-independent analytic
  set is contained in a $G_{i,j}$-independent \Borel set, Lemmas \ref
  {graph:dichotomy:splitting:one} and \ref
  {graph:dichotomy:splitting:two} ensure that if $\sequence{X_{i,j}}[i,j
  \in \N]$ is a sequence of analytic sets and $a$ is an $\sequence
  {X_{i.j}}[i,j \in \N]$-terminal approximation, then there is a \Borel set $B(a,
  \sequence{X_{i,j}}[i,j \in \N]) \supseteq A(a, \sequence{X_{i,j}}[i,j \in
  \N])$ that is $G_{k_{n^a}, f(k_{n^a})}$-independent if $k_{n^a} \in
  \domain{f^a}$, and $G_{k_{n^a}, \ell}$-independent for all $\ell \in
  \N$ if $k_{n^a} \notin \domain{f^a}$.
  
  We finally define $X_{k,\ell}^{\alpha + 1}$ to be the difference of
  $X_{k,\ell}^\alpha$ and the union of the sets of the form $B(a,
  \sequence{X_{i,j}^\alpha}[i,j \in \N])$, where $a$ is an $\sequence
  {X_{i,j}^\alpha}[i,j \in \N]$-terminal approximation, $k_{n^a} = k$,
  and $f^a(k_{n^a}) \ge \ell$ if $k_{n^a} \in \domain{f^a}$.
  
  \begin{lemma} \label{graph:dichotomy:extension}
    Suppose that $\alpha < \omega_1$ and $a$ is an approximation
    that is not $\sequence{X_{i,j}^{\alpha+1}}[i,j \in \N]$-terminal. Then
    there is a one-step extension of $a$ that is not $\sequence
    {X_{i,j}^{\alpha}}[i,j \in \N]$-terminal.
  \end{lemma}
  
  \begin{lemmaproof}
    Fix a one-step extension $b$ of $a$ for which there is a
    configuration $\gamma$ compatible with $b$ and $\sequence
    {X_{i,j}^{\alpha+1}}[i,j \in \N]$. Note that if $k_{n^b} \in \domain
    {f^b}$, then $(\phi_X \composition \phi^\gamma)(s_{n^b}) \in
    X_{k_{n^b}, f^b(k_{n^b})}^{\alpha+1}$, so $A(b, \sequence{X_{i,j}^
    \alpha}[i,j \in \N]) \intersection X_{k_{n^b}, f^b(k_{n^b})}^{\alpha+1}
    \neq \emptyset$, thus $b$ is not $\sequence{X_{i,j}^\alpha}[i,j \in
    \N]$-terminal. And if $k_{n^b} \notin \domain{f^b}$, then there
    exists $\ell \in \N$ for which $(\phi_X \composition \phi^\gamma)
    (s_{n^b}) \in X_{k_{n^b}, \ell}^{\alpha + 1}$, so $A(b, \sequence
    {X_{i,j}^\alpha}[i,j \in \N]) \intersection X_{k_{n^b}, \ell}^{\alpha+1}
    \neq \emptyset$, thus $b$ is not $\sequence{X_{i,j}^\alpha}[i,j \in
    \N]$-terminal.
  \end{lemmaproof}
  
  Fix $\alpha < \omega_1$ such that the families of $\sequence
  {X_{i,j}^\alpha}[i,j \in \N]$-terminal approximations and $\sequence
  {X_{i,j}^{\alpha+1}}[i,j \in \N]$-terminal approximations are the
  same, let $a_0$ denote the unique approximation $a$ with the
  property that $n^a = 0$, and observe that $A(a_0, \sequence
  {X_{i,j}}[i, j \in \N]) = \intersection[i \in \N][{\union[j \in \N][X_{i,j}]}]$
  for all sequences $\sequence{X_{i,j}}[i,j \in \N]$ of subsets of $X$.
  In particular, it follows that if $a_0$ is $\sequence{X_{i,j}^\alpha}[i,j
  \in \N]$-terminal, then $\intersection[i \in \N][{\union[j \in \N]
  [X^{\alpha+1}_{i,j}]}] = \emptyset$, so condition (1) holds.
  
  Otherwise, by recursively applying Lemma \ref
  {graph:dichotomy:extension}, we obtain one-step extensions
  $a_{n+1}$ of $a_n$ that are not $\sequence{X_{i,j}^\alpha}[i,j \in
  \N]$-terminal for all $n \in \N$. Define $f \from \N \to \N$ by $f =
  \union[n \in \N][f^{a_n}]$, define $\phi
  \from \Cantorspace \to \Bairespace$ by $\phi(c) = \union[n \in \N]
  [\phi^{a_n}(\restriction{c}{n})]$ for all $c \in \Cantorspace$, and
  define $\psi_n \from \Cantorspace \to \Bairespace$ by $\psi_n(c)
  = \union[m \in \N][\psi_n^{a_{n+1+m}}(\restriction{c}{m})]$ for all $c
  \in \Cantorspace$ and $n \in \N$. To see that $\phi_X \composition
  \phi$ is a homomorphism from $\sequence{\Gzero[k]}[k \in \N]$ to
  $\sequence{G_{k, f(k)}}[k \in \N]$, we will show that $(\phi_{k_n,
  f(k_n)} \composition \psi_n)(c) = \sequence{(\phi_X \composition
  \phi)(s_n \concatenation \sequence{i} \concatenation c)}[i < 2]$ for
  all $c \in \Cantorspace$ and $n \in \N$. For this, it is sufficient to show that if $U \subseteq X \times X$ is an open neighborhood of $(\phi_{k_n, f(k_n)}
  \composition \psi_n)(c)$ and $V \subseteq X \times X$ is an open neighborhood of
  $\sequence{(\phi_X \composition \phi)(s_n \concatenation
  \sequence{i} \concatenation c)}[i < 2]$, then $U \intersection V \neq
  \emptyset$. Towards this end, fix
  $m \in \N$ for which $\image
  {\phi_{k_n, f(k_n)}}{\extensions{\psi_n^{a_{n+1+m}}(s)}} \subseteq
  U$ and $\product[i <
  2][\image{\phi_X}{\extensions{\phi^{a_{n+1+m}}(s_n \concatenation
  \sequence{i} \concatenation s)}}] \subseteq V$, where $s = \restriction{c}{m}$. The fact that $a_m$ is not
  $\sequence{X_{i,j}^\alpha}[i,j \in \N]$-terminal then yields a
  configuration $\gamma$ compatible with $a_m$, so $(\phi_{k_n, f(k_n)}
  \composition \psi_n^\gamma)(s) \in U$ and $\sequence
  {(\phi_X \composition \phi^\gamma)(s_n \concatenation \sequence
  {i} \concatenation s)}[i < 2] \in V$, thus $U \intersection V
  \neq \emptyset$.
\end{theoremproof}

\section{Lacunary sets} \label{lacunary}

Here we note the connection between condition (1) of Theorem \ref
{graph:dichotomy} and lacunary sets.

\begin{proposition} \label{lacunary:lacunarytoindependent}
  Suppose that $\Gamma$ is a tsi analytic \Hausdorff group, $X$ is an analytic
  \Hausdorff space, $\Gamma \action X$ is a $\sigma$-lacunary
  \Borel action such that $\orbitrelation{\Delta}{X}$ is \Borel for all
  open sets $\Delta \subseteq \Gamma$, $\sequence{\Delta_i}[i \in
  \N]$ is a neighborhood basis of $\identity{\Gamma}$ consisting of
  conjugation-invariant symmetric open sets, and $G_{i,j} =
  \orbitrelation{\Delta_i}{X} \setminus \orbitrelation{\Delta_j}{X}$ for
  all $i, j \in \N$. Then there are \Borel sets $B_n \subseteq X$ such
  that $X = \union[n \in \N][B_n]$ and $\forall n \in \N \exists i \in \N
  \forall j \in \N \ \Borelchromaticnumber{\restriction{G_{i,j}}{B_n}} \le
  \aleph_0$.
\end{proposition}

\begin{propositionproof}
  By breaking $X$ into countably-many $\orbitequivalencerelation
  {\Gamma}{X}$-invariant \Borel sets, we can assume that there is
  an open neighborhood $\Delta \subseteq \Gamma$ of $\identity
  {\Gamma}$ for which there is a $\Delta$-lacunary
  $\orbitequivalencerelation{\Gamma}{X}$-complete \Borel set $B
  \subseteq X$.
  
  Fix $i \in \N$ for which there is an open neighborhood $\Delta'
  \subseteq \Gamma$ of $\identity{\Gamma}$ such that $\inverse
  {(\Delta')} \Delta_i \Delta' \subseteq \Delta$. To see that
  $\Borelchromaticnumber{G_{i,j}} \le \aleph_0$ for all $j \in \N$,
  fix $j \in \N$ and an open set $\Delta'' \subseteq \Delta'$ such
  that $\Delta'' \inverse{(\Delta'')} \subseteq \Delta_j$.
  
  \begin{lemma}
    The set $\Delta'' B$ is $G_{i,j}$-independent.
  \end{lemma}
  
  \begin{lemmaproof}
    Suppose that $x'', y'' \in \Delta'' B$ are $\orbitrelation{\Delta_i}
    {X}$-related. Then there exist $\delta_x'', \delta_y'' \in \Delta''$ for
    which the points $x = \inverse{(\delta_x'')} \cdot x''$ and $y =
    \inverse{(\delta_y'')} \cdot y''$ are in $B$. As $x$ and $y$ are
    $\orbitrelation{\inverse{(\Delta'')} \Delta_i \Delta''}{X}$-related, so
    $\orbitrelation{\Delta}{X}$-related, thus equal, it follows that $x''$
    and $y''$ are $\orbitrelation{\Delta'' \inverse{(\Delta'')}}{X}$-related,
    thus $\orbitrelation{\Delta_j}{X}$-related.
  \end{lemmaproof}
  
  The conjugation invariance of $\Delta_i$ and $\Delta_j$ now
  ensures that $\gamma \Delta'' B$ is $G_{i,j}$-independent, and
  therefore contained in an $G_{i,j}$-independent \Borel set, for all
  $\gamma \in \Gamma$. As $X$ is the union of countably-many
  sets of this form, it follows that $\Borelchromaticnumber{G_{i,j}}
  \le \aleph_0$.
\end{propositionproof}

A topological group is \definedterm{cli} if it has a compatible
complete left-invariant metric, or equivalently, a compatible
complete right-invariant metric (see, for example, \cite[Proposition
3.A.2]{Becker}). It is well-known that every tsi group is cli (see, for
example, \cite[Corollary 1.2.2]{BeckerKechris}).

\begin{proposition} \label{lacunary:independenttolacunary}
  Suppose that $\Gamma$ is a cli \Polish group, $X$ is an analytic
  metric space, $\Gamma \action X$ is continuous, $\sequence
  {\Delta_i}[i \in \N]$ is a neighborhood basis of $\identity{\Gamma}$
  consisting of symmetric open sets, $G_{i,j} = \orbitrelation
  {\Delta_i}{X} \setminus \orbitrelation{\Delta_j}{X}$ for all $i, j
  \in \N$, and there are \Borel sets $B_n \subseteq X$ with the
  property that $X = \union[n \in \N][B_n]$ and $\forall n \in \N \exists
  i \in \N \forall j \in \N \ \Borelchromaticnumber{\restriction{G_{i,j}}
  {B_n}} \le \aleph_0$. Then $\Gamma \action X$ is
  $\sigma$-lacunary.
\end{proposition}

\begin{propositionproof}
  We can assume that $\Gamma$ is not discrete, since otherwise
  $\Gamma \action X$ is trivially $\sigma$-lacunary. So by passing
  to a subsequence of $\sequence{\Delta_i}[i \in \N]$, we can also
  assume that $\closure{\Delta_{i+1}}^2 \subseteq \Delta_i$ for all
  $i \in \N$. By breaking each $B_n$ into countably-many \Borel
  sets, we can moreover assume that there are natural numbers
  $i_n \in \N$ such that $B_n$ is $G_{i_n, i_n + 3}$-independent and
  $\Borelchromaticnumber{\restriction{G_{i_n, i_n + 4 + j}}{B_n}} \le
  \aleph_0$ for all $j, n \in \N$. As a result of \Montgomery-\Novikov
  ensures that the class of \Borel sets is closed under category
  quantification (see, for example, \cite[Theorem 16.1]{Kechris}), it
  follows that the map $\phi \from X \to \N$ given by $\phi(x) =
  \min \set{n \in \N}[\existnonmeagerlymany \gamma \in \Gamma
  \ \gamma \cdot x \in B_n]$ is \Borel. By passing to the
  $\orbitequivalencerelation{\Gamma}{X}$-invariant \Borel sets $X_n
  = \preimage{\phi}{B_n}$, it is sufficient to show that if $i \in \N$ and
  there is a $G_{i,i+3}$-independent \Borel set $B \subseteq X$ with
  the property that $\forall j \in \N \ \Borelchromaticnumber{\restriction
  {G_{i,i+4+j}}{B}} \le \aleph_0$ and $\forall x \in X
  \existnonmeagerlymany \gamma \in \Gamma \ \gamma \cdot x \in
  B$, then there is a $\Delta_{i+2}$-lacunary
  $\orbitequivalencerelation{\Gamma}{X}$-complete \Borel set.
  
  Towards this end, observe that the set $E = \restriction{\orbitrelation
  {\Delta_{i+3}}{X}}{B}$ is an equivalence relation. As $E$ has
  countable index below $\restriction{\orbitequivalencerelation
  {\Gamma}{X}}{B}$, by thinning down $B$ if necessary, we can
  assume that $\forall x \in B \existnonmeagerlymany \gamma \in
  \Gamma \ x \mathrel{E} \gamma \cdot x$. Fix positive real numbers
  $\epsilon_j \goesto 0$, as well as \Borel colorings $c_{i+4+j} \from
  B \to \N$ of $\restriction{G_{i,i+4+j}}{B}$ such that $\diameter
  {\preimage{c_{i+4+j}}{\set{m}}} \le \epsilon_j$ for all $j, m \in \N$.
  For each $j \in \N$ and $x \in B$, let $s_{i+4+j}(x)$ denote the
  lexicographically minimal sequence $s \in \Bairespace[j+1]$ for
  which there are non-meagerly many $\gamma \in \Gamma$ with
  the property that $\gamma \cdot x \in \intersection[k \le j][\preimage
  {c_{i+4+k}}{\set{s(k)}}] \intersection \equivalenceclass{x}{E}$, and
  let $C_{i+4+j}$ denote the set of $x \in B$ for which $s_{i+4+j}(x) =
  \sequence{c_{i+4+k}(x)}[k \le j]$.
  
  A \definedterm{ray} from $x \in B$ through $\sequence{C_{i+4+j}}[j
  \in \N]$ is a sequence $\sequence{\delta_{i+3+j}}[j \in \N]$ with the
  property that $\delta_{i+3+j} \in \Delta_{i+3+j}$ and $\delta_{i+3+j}
  \cdots \delta_{i+3} \cdot x \in C_{i+4+j}$ for all $j \in \N$. A
  straightforward recursive construction yields the existence of such
  rays, while a straightforward inductive argument ensures that if
  $\sequence{\delta_{i+3+j}}[j \in \N]$ is such a ray, then $\delta_{i+
  3+k} \cdots \delta_{i+3+j} \in \Delta_{i+2+j}$ for all $k > j$. In
  particular, it follows that $\sequence{\delta_{i+3+j} \cdots \delta_{i+
  3}}[j \in \N]$ is \Cauchy with respect to every compatible complete
  right-invariant metric on $\Gamma$, and therefore converges to
  some $\delta \in \closure{\Delta_{i+2}}$.

  Observe now that if $\sequence{\delta^x_{i+3+j}}[j \in \N]$ and
  $\sequence{\delta^y_{i+3+j}}[j \in \N]$ are rays from points $x$ and
  $y$ in $B$ through $\sequence{C_{i+4+j}}[j \in \N]$, and $\delta^x$
  and $\delta^y$ are the corresponding limit points, then $\delta^x
  \cdot x \mathrel{\orbitrelation{\Delta_{i+2}}{X}} \delta^y \cdot y
  \implies x \mathrel{\orbitrelation{\Delta_i}{X}} y \implies x \mathrel{E}
  y$ and $x \mathrel{E} y \implies \delta^x \cdot x = \delta^y \cdot y$.
  We therefore obtain a function $\psi \from B \to X$ by insisting that
  $\psi(x) = y$ if and only if there is a ray $\sequence{\delta_{i+3+j}}[j
  \in \N]$ from $x$ through $\sequence{C_{i+4+j}}[j \in \N]$ for which
  $\delta_{i+3+j} \cdots \delta_{i+3} \cdot x \goesto y$. It also follows
  that the corresponding set $\image{\psi}{B}$ is $\Delta_{i+
  2}$-lacunary, and the fact that $\forall y \in \image{\psi}{B}
  \existnonmeagerlymany \gamma \in \Gamma \ \psi(\gamma \cdot y)
  = y$ ensures that $\image{\psi}{B}$ is \Borel.
\end{propositionproof}

\section{Compositions} \label{compositions}

Here we note several ways of refining condition (2) of Theorem \ref
{graph:dichotomy}.

\begin{proposition} \label{compositions:alignment}
  Suppose that $f \from \N \to \N$. Then there is a continuous
  homomorphism $\phi \from \Cantorspace \to \Cantorspace$
  from $\sequence{\Gzero[k]}[k \in \N]$ to $\sequence{\Gzero[f(k)]}
  [k \in \N]$.
\end{proposition}

\begin{propositionproof}
  Recursively construct $m_n \in \N$ and $u_n \in \Cantortree$ with
  the property that $k_{m_n} = f(k_n)$ and $s_{m_n} = \phi_n(s_n)$,
  where $\phi_n \from \Cantorspace[n] \to \Cantorspace[m_n]$ is
  given by $\phi_n(t) = u_0 \concatenation \concatenation[i < n][t(i)
  \concatenation u_{i+1}]$ for all $t \in \Cantorspace[n]$, and define
  $\phi \from \Cantorspace \to \Cantorspace$ by $\phi(c) = \union[n
  \in \N][\phi_n(\restriction{c}{n})]$ for all $c \in \Cantorspace$.
  
  To see that $\phi$ is a homomorphism from $\sequence{\Gzero[k]}
  [k \in \N]$ to $\sequence{\Gzero[f(k)]}[k \in \N]$, observe that if $c
  \in \Cantorspace$ and $n \in \N$, then there exists $d \in
  \Cantorspace$ such that $\phi(s_n \concatenation \sequence{i}
  \concatenation c) = s_{m_n} \concatenation \sequence{i}
  \concatenation d$ for all $i < 2$. As $k_{m_n} = f(k_n)$, it follows
  that $\phi(s_n \concatenation \sequence{0} \concatenation c)
  \mathrel{\Gzero[f(k_n)]} \phi(s_n \concatenation \sequence{1}
  \concatenation c)$.
\end{propositionproof}

For all $s, t \in \Cantortree$, we use $\involutiongraph{s}[t]$ to
denote the subgraph of $\involutiongraph{s}$ given by
$\involutiongraph{s}[t] = \set{\sequence{s \concatenation \sequence
{i} \concatenation t \concatenation c}[i < 2]}[c \in \Cantorspace]$.

\begin{proposition} \label{compositions:focus}
  Suppose that $\sequence{R_{i,j}}[i,j \in \N]$ is a sequence of
  analytic binary relations on $\Cantorspace$ with the property that
  $\Gzero[k] \subseteq \union[j \in \N][R_{i,j}]$ for all $i, k \in \N$.
  Then there are functions $g_n \from \Cantorspace[<n] \to \N$ and
  a continuous homomorphism $\phi \from \Cantorspace \to
  \Cantorspace$ from $\sequence{\Gzero[k]}[k \in \N]$ to $\sequence
  {\Gzero[k]}[k \in \N]$ that is also a homomorphism from $\sequence
  {\involutiongraph{s_n}[t]}[n \in \N, t \in \Cantortree]$ to $\sequence
  {R_{k_{n+1+\length{t}}, g_{n+1+\length{t}}(t)}}[n \in \N, t \in
  \Cantortree]$.
\end{proposition}

\begin{propositionproof}
  We will recursively construct $g_n \from \Cantortree[n] \to \N$,
  $m_n \in \N$, $u_n \in \Cantortree$, and open sets $U_{j,n}
  \subseteq \Cantorspace$, from which we define
  $\phi_{\closedopeninterval{m}{n}} \from \Cantorspace[n-m] \to
  \Cantortree$ by $\phi_{\closedopeninterval{m}{n}}(t) =
  \concatenation[i < n - m][u_{i+m} \concatenation \sequence{t(i)}]$
  for all $m \le n$ and $t \in \Cantorspace[n-m]$, satisfying the
  following conditions:
  \begin{enumerate}
    \item $\forall n \in \N \forall t \in \Cantortree[n] \forall c \in
      \intersection[j \in \N][U_{j,n}] \\
        \hspace*{5pt} \sequence{\phi_{\closedopeninterval{0}{n}}(s_{n-
          1 - \length{t}} \concatenation \sequence{i} \concatenation t)
            \concatenation c}[i < 2] \in R_{k_n, g_n(t)}$.
    \item $\forall j, n \in \N \ U_{j, n}$ is dense in $\extensions{u_n}$.
    \item $\forall n \in \N \forall j \le n \forall t \in \Cantorspace[j+1]
      \ \extensions{\phi_{\closedopeninterval{n - j}{n + 1}}(t)} \subseteq
        U_{j, n - j}$.
    \item $\forall n \in \N \ (k_{m_n} = k_n \mathand s_{m_n} =
      \phi_{\closedopeninterval{0}{n}} (s_n) \concatenation u_n)$.
  \end{enumerate}
  
  Suppose that $n \in \N$ and we have already found $g_k$, $m_k$,
  $u_k$, and $\sequence{U_{j,k}}[j \in \N]$ for all $k < n$. For all $g
  \from \Cantorspace[<n] \to \N$, let $B_g$ be the set of $c \in
  \Cantorspace$ such that $\sequence{\phi_{\closedopeninterval{0}
  {n}}(s_{n-1-\length{t}} \concatenation \sequence{i} \concatenation
  t) \concatenation c}[i < 2] \in R_{k_n, g(t)}$ for all $t \in
  \Cantortree[n]$. Fix $g_n \from \Cantorspace[<n] \to \N$ for which
  $B_{g_n}$ is non-meager, as well as $u_{0,n} \in \Cantortree$ for
  which $B_{g_n}$ is comeager in $\extensions{u_{0,n}}$, in addition
  to dense open sets $U_{j,n} \subseteq \extensions{u_{0,n}}$ for
  which $\intersection[j \in \N][U_{j,n}] \subseteq B_{g_n}$. Fix an
  enumeration $\sequence{v_{k,n}}[k < \ell]$ of $\Cantorspace[\le
  n]$, and recursively find extensions $u_{k+1,n} \in \Cantortree$ of
  $u_{k, n}$ such that $\extensions{\phi_{\closedopeninterval{n -
  \length{v_{k,n}}}{n}}(v_{k,n}) \concatenation u_{k+1,n}} \subseteq
  U_{\length{v_{k,n}},n-\length{v_{k,n}}}$ for all $k < \ell$. Finally, fix
  $m_n \in \N$ and an extension $u_n \in \Cantortree$ of $u_{\ell,n}$
  for which $k_{m_n} = k_n$ and $s_{m_n} = \phi_{\closedopeninterval
  {0}{n}}(s_n) \concatenation u_n$.
  
  Define $\phi_{\closedopeninterval{m}{\infty}} \from \Cantorspace \to
  \Cantorspace$ by $\phi_{\closedopeninterval{m}{\infty}}(c) = \union
  [n \in \N][\phi_{\closedopeninterval{m}{m+n}}(\restriction{c}{n})]$ for
  all $c \in \Cantorspace$ and $m \in \N$. Condition (4) ensures that
  $\phi_{\closedopeninterval{0}{\infty}}$ is a homomorphism from
  $\sequence{\Gzero[k]}[k \in \N]$ to $\sequence{\Gzero[k]}[k \in \N]$.
  To see that $\phi_{\closedopeninterval{0}{\infty}}$ is a
  homomorphism from $\sequence{\involutiongraph{s_n}[t]}[n \in \N,
  t \in \Cantortree]$ to $\sequence{R_{k_{n+1+\length{t}}, g_{n+1+
  \length{t}}(t)}}[n \in \N, t \in \Cantortree]$, suppose that $c \in
  \Cantorspace$, $n \in \N$, and $t \in \Cantortree$. Condition (3)
  then ensures that $\extensions{\phi_{\closedopeninterval{n+1+
  \length{t}}{n+1+\length{t}+j+1}}(\restriction{c}{(j+1)})} \subseteq
  U_{j,n+1+\length{t}}$ for all $j \in \N$, so
  $\phi_{\closedopeninterval{n+1+\length{t}}{\infty}}(c) \in
  \intersection[j \in \N][U_{j,n+1+\length{t}}]$, in which case condition
  (1) implies that $\sequence{\phi_{\closedopeninterval{0}{n+1+
  \length{t}}}(s_n \concatenation \sequence{i} \concatenation t)
  \concatenation \phi_{\closedopeninterval{n+1+\length{t}}{\infty}}
  (c)}[i < 2] \in R_{k_{n+1+\length{t}}, g_{n+1+\length{t}}(t)}$. But
  $\phi_{\closedopeninterval{0}{\infty}}(s_n \concatenation \sequence
  {i} \concatenation t \concatenation c) = \phi_{\closedopeninterval{0}
  {n+1+\length{t}}}(s_n \concatenation \sequence{i} \concatenation t)
  \concatenation \phi_{\closedopeninterval{n+1+\length{t}}{\infty}}(c)$
  for all $i < 2$, thus $\sequence{\phi_{\closedopeninterval{0}{\infty}}
  (s_n \concatenation \sequence{i} \concatenation t \concatenation c)}
  [i < 2] \in R_{k_{n+1+\length{t}}, g_{n+1+\length{t}}(t)}$.
\end{propositionproof}

For all $F \subseteq \N \times \N$ and $c, d \in \Cantorspace[F]$,
let $\differenceset{c}{d}$ be the set of $\pair{m}{n} \in F$ with
$c(m, n) \neq d(m, n)$. For all $i \in \N$, set $\differenceset[i]{c}
{d} = \differenceset{c}{d} \intersection (i \times \N)$. When $F \in
\sets{<\aleph_0}{i \times \N}$, set $\D{i}[F] = \set{\pair{c}{d} \in
\Cantorspace[\N \times \N] \times \Cantorspace[\N \times \N]}
[{\differenceset[i]{c}{d} = F}]$.

\begin{proposition} \label{compositions:meager}
  Suppose that $i \in \N$, $F \in \sets{<\aleph_0}{i \times \N}$, $R
  \subseteq \D{i}[F]$ has the \Baire property, and there are densely
  many $u \in \Cantorspace[<(\N \times \N)]$ for which there is a
  homeomorphism $\phi \from \extensions{u} \to \extensions{u}$
  whose graph is contained in $\D{i}[\emptyset] \setminus R
  \inverse{R}$. Then $R$ is meager.
\end{proposition}

\begin{propositionproof}
  Suppose, towards a contradiction, that $R$ is non-meager. Then
  there exist $G \in \sets{<\aleph_0}{(i \times \N) \setminus F}$ and
  $H, H' \in \sets{<\aleph_0}{(\N \setminus i) \times \N}$ for which
  there exist $r \in \Cantorspace[F]$, $s \in \Cantorspace[G]$, $t \in
  \Cantorspace[H]$, and $t' \in \Cantorspace[H']$ with the property
  that $R$ is comeager in $\D{i}[F] \intersection (\extensions{r \union
  s \union t} \times \extensions{\flip{r} \union s \union t'})$, in which
  case the set $S$ of $\pair{c}{\pair{d}{d'}} \in \Cantorspace[(i \times
  \N) \setminus (F \union G)] \times (\Cantorspace[((\N \setminus i)
  \times \N) \setminus H] \times \Cantorspace[((\N \setminus i) \times
  \N) \setminus H'])$ with the property that $((c \union r \union s)
  \union (d \union t)) \mathrel{R} ((c \union \flip{r} \union s) \union (d'
  \union t'))$ is comeager. 
  
  Let $C$ denote the set of $c \in \Cantorspace[(i \times \N)
  \setminus (F \union G)]$ for which $\verticalsection{S}{c}$ is
  comeager, and let $D$ denote the set of $\pair{c}{d} \in
  \Cantorspace[(i \times \N) \setminus (F \union G)] \times
  \Cantorspace[((\N \setminus i) \times \N) \setminus H]$ for which
  $\verticalsection{(\verticalsection{S}{c})}{d}$ is comeager. The
  \Kuratowski-\Ulam theorem ensures that $C$ is comeager, as is
  $\verticalsection{D}{c}$ for all $c \in C$, thus $\verticalsection{D}{c}
  \times \verticalsection{D}{c} \subseteq \verticalsection{S}{c}
  \inverse{\verticalsection{S}{c}}$ for all $c \in C$. 
  
  Fix $I \in \sets{<\aleph_0}{(i \times \N) \setminus(F \union G)}$
  and $J \in \sets{<\aleph_0}{((\N \setminus i) \times \N) \setminus
  H}$ for which there exist $u \in \Cantorspace[I]$ and $v \in
  \Cantorspace[J]$ with the property that there is a homeomorphism
  $\phi \from \extensions{(r \union s \union u) \union (t \union v)} \to
  \extensions{(r \union s \union u) \union (t \union v)}$ whose graph is
  contained in $\D{i, \emptyset} \setminus R \inverse{R}$. Fix $c \in
  C \intersection\extensions{u}$ and define $\psi \from \Cantorspace
  [((\N \setminus i) \times \N) \setminus H] \intersection \extensions
  {v} \to \Cantorspace[((\N \setminus i) \times \N) \setminus H]
  \intersection \extensions{v}$ by $\psi(d) = (\projection
  [{\Cantorspace[((\N \setminus i) \times \N) \setminus H]}]
  \composition \phi)((c \union r \union s) \union (d \union t))$ for all $d
  \in \Cantorspace[((\N \setminus i) \times \N) \setminus H]
  \intersection \extensions{v}$. The fact that $\psi$ is a
  homeomorphism then ensures that there are comeagerly many
  $d \in \Cantorspace[((\N \setminus i) \times \N) \setminus H]
  \intersection \extensions{v}$ that are also in $\verticalsection{D}{c}
  \intersection \preimage{\psi}{\verticalsection{D}{c}}$. But the
  defining property of $\phi$ ensures that $d$ and $\psi(d)$ are not
  $(\verticalsection{S}{c} \inverse{\verticalsection{S}{c}})$-related,
  the desired contradiction.
\end{propositionproof}

For each $i \in \N$, define $\differencefunction{i} \from
\CantorCantorspace \times \CantorCantorspace \to \N \union \set
{\aleph_0}$ by setting $\differencefunction{i}(c, d) = \cardinality
{\differenceset{c}{d} \intersection (\set{i} \times \N)}$ for all $c, d \in
\CantorCantorspace$.

A \definedterm{homomorphism} from a function $f \from X \times X
\to N$ to a function $g \from Y \times Y \to N$ is a map $\phi \from
X \to Y$ such that $f(w, x) = g(\phi(w), \phi(x))$ for all $w, x \in X$.
More generally, a \definedterm{homomorphism} from a sequence
$\sequence{f_i \from X \times X \to N}[i \in I]$ to a sequence
$\sequence{g_i \from Y \times Y \to N}[i \in I]$ is a map $\phi \from
X \to Y$ that is a homomorphism from $f_i$ to $g_i$ for all $i \in I$.

\begin{proposition} \label{compositions:continuous}
  Suppose that $C \subseteq \CantorCantorspace$ is comeager.
  Then there is a continuous homomorphism $\phi \from
  \CantorCantorspace \to C$ from $\sequence{\delta_i}[i \in \N]$
  to $\sequence{\delta_i}[i \in \N]$.
\end{proposition}

\begin{propositionproof}
  Fix dense open sets $U_n \subseteq \CantorCantorspace$ for
  which $\intersection[n \in \N][U_n] \subseteq C$.
  
  \begin{lemma} \label{compositions:continuous:innerrecursion}
    For all $F, G \in \sets{<\aleph_0}{\N \times \N}$, $\phi \from
    \Cantorspace[F] \to \Cantorspace[G]$, and $n \in \N$, there
    exist $H \in \sets{<\aleph_0}{\protect \setcomplement{G}}$ and
    $t \in \Cantorspace[H]$ such that $\extensions
    {\phi(s) \union t} \subseteq U_n$ for all $s \in \Cantorspace[F]$.
  \end{lemma}
  
  \begin{lemmaproof}
    Fix an enumeration $\sequence{s_m}[m < 2^{\cardinality{F}}]$
    of $\Cantorspace[F]$, and recursively find pairwise disjoint sets
    $H_m \in \sets{<\aleph_0}{\protect \setcomplement{G}}$ and $t_m
    \in \Cantorspace[H_m]$ with $\extensions{\phi(s_m) \union
    \union[\ell \le m][t_\ell]} \subseteq U_n$ for all $m < 2^{\cardinality{F}}$.
    Define $H = \union[m < 2^{\cardinality{F}}][H_m]$ and $t = \union
    [m < 2^{\cardinality{F}}][t_m]$.
  \end{lemmaproof}
   
  Fix an injective enumeration $\sequence{i_n, j_n}[n \in \N]$ of
  $\N \times \N$, and for all $n \in \N$, set $F_n = \set{\pair{i_m}
  {j_m}}[m < n]$. By recursively appealing to Lemma \ref
  {compositions:continuous:innerrecursion}, we obtain $H_n \in
  \sets{<\aleph_0}{\N \times \N}$ and $j_n' \in \setcomplement
  {\verticalsection{(H_n)}{i_n}}$ for which the sets $G_n = H_n \union
  \set{\pair{i_n}{j_n'}}$ are pairwise disjoint, as well as $t_n \in
  \Cantorspace[H_n]$ such that $\extensions{\phi_n(s) \union t_n}
  \subseteq U_n$ for all $n \in \N$ and $s \in \Cantorspace[F]$, where
  $\phi_n \from \Cantorspace[F_n] \to \Cantorspace[{\union[m < n]
  [G_m]}]$ is given by $\phi_n(s) = \union[m < n][t_{s(i_m,j_m),m}]$,
  and $t_{k,m}$ is the extension of $t_m$ sending $\pair{i_n}{j_n'}$ to
  $k$, for all $k < 2$ and $m \in \N$. Then the function $\phi \from
  \Cantorspace[\N \times \N] \to \Cantorspace[\N \times \N]$, obtained
  by insisting that $\support{\phi(c)} \subseteq \union[n \in \N][G_n]$
  and $\restriction{\phi(c)}{G_n} = t_{c(i_n, j_n), n}$ for all $c \in
  \CantorCantorspace$ and $n \in \N$, is continuous.
  
  To see that $\phi$ is a homomorphism from $\delta_i$ to $\delta_i$
  for all $i \in \N$, simply observe that $\differenceset{\phi(c)}{\phi(d)}
  = \set{\pair{i_n}{j_n'}}[n \in \N \mathand \pair{i_n}{j_n} \in
  \differenceset{c}{d}]$ for all $c, d \in \CantorCantorspace$.
\end{propositionproof}

Let $\equality{X}$ denote equality on $X$. We will abuse notation by
identifying $\Ethree$, and more generally $\functions{k}{\equality
{\Cantorspace}} \times \Ethree$ for all $k \in \N$, with the corresponding
equivalence relations on $\CantorCantorspace$.

\begin{proposition} \label{compositions:embedding}
  Suppose that $D \subseteq \CantorCantorspace \times
  \CantorCantorspace$ is closed and nowhere dense in $\D{i}[F]$
  for all $i \in \N$ and $F \in \sets{<\aleph_0}{i \times \N}$, and $R
  \subseteq \CantorCantorspace \times \CantorCantorspace$ is
  meager in $\D{i}[F]$ for all $i \in \N$ and $F \in \sets{<\aleph_0}{i
  \times \N}$. Then there is a continuous homomorphism $\phi \from
  \Cantorspace[\N \times \N] \to \Cantorspace[\N \times \N]$ from
  $\sequence{\equality{\Cantorspace}^k \times \Ethree}[k \in \N]$ to
  $\sequence{\equality{\Cantorspace}^k \times \Ethree}[k \in \N]$ that
  is also a homomorphism from $\pair{\setcomplement{\equality
  {\CantorCantorspace}}}{\setcomplement{\Ethree}}$ to $\pair
  {\setcomplement{D}}{\setcomplement{R}}$.
\end{proposition}

\begin{propositionproof}
  For all $i \in \N$ and $F \in \sets{<\aleph_0}{i \times \N}$, fix a
  decreasing sequence $\sequence{U_{i, F, n}}[n \in \N]$ of
  dense open symmetric subsets of $\D{i}[F] \setminus D$ whose
  intersection is disjoint from $R$.
  
  \begin{lemma} \label{compositions:embedding:innerrecursion}
    For all $F, G \in \sets{<\aleph_0}{\N \times \N}$, $\phi \from
    \Cantorspace[F] \to \Cantorspace[G]$, and $i, n \in \N$, there
    exist $H \in \sets{<\aleph_0}{\protect \setcomplement{G}}$ and
    $t_0, t_1 \in \Cantorspace[H]$ with the property that
    $\differenceset[i]{t_0}{t_1} = \emptyset$ and $\D{i}[{\differenceset
    [i]{\phi(s_0)}{\phi(s_1)}}] \intersection \product[k < 2][\extensions
    {\phi(s_k) \union t_k}] \subseteq U_{i, \differenceset[i]{\phi(s_0)}
    {\phi(s_1)}, n}$ for all $s_0, s_1 \in \Cantorspace[F]$.
  \end{lemma}
  
  \begin{lemmaproof}
    Fix an enumeration $\sequence{s_{0,m}, s_{1,m}}[m < 4^
    {\cardinality{F}}]$ of $\Cantorspace[F] \times \Cantorspace[F]$,
    and recursively find pairwise disjoint sets $H_m \in \sets{<
    \aleph_0}{\protect \setcomplement{G}}$ and $t_{0,m}, t_{1,m}
    \in \Cantorspace[H_m]$ such that $\differenceset[i]{t_{0,m}}{t_{1,m}}
    = \emptyset$ and
    \begin{equation*}
      \textstyle
      \D{i}[{\differenceset[i]{\phi(s_{0,m})}{\phi(s_{1,m})}}] \intersection
        \product[k < 2][\extensions{\phi(s_{k,m}) \union \union[\ell \le
          m][t_{k,\ell}]}] \subseteq U_{i, \differenceset[i]{\phi(s_{0,m})}
            {\phi(s_{1,m})}, n}
    \end{equation*}
    for all $m < 4^{\cardinality{F}}$. Set $H = \union[m < 4^{\cardinality
    {F}}][H_m]$ and $t_k = \union[m < 4^{\cardinality{F}}][t_{k, m}]$.
  \end{lemmaproof}
  
  Fix an injective enumeration $\sequence{i_n, j_n}[n \in \N]$ of
  $\N \times \N$, and for all $n \in \N$, set $F_n = \set{\pair{i_m}
  {j_m}}[m < n]$. By recursively appealing to Lemma \ref
  {compositions:embedding:innerrecursion}, we obtain pairwise
  disjoint sets $G_n \in \sets{<\aleph_0}{\N \times \N}$ and
  $t_{0, n}, t_{1, n} \in \Cantorspace[G_n]$ such that
  $\differenceset[i_n]{t_{0,n}}{t_{1,n}} = \emptyset$ and
  \begin{equation*}
    \textstyle
    \D{i_n}[{\differenceset[i_n]{\phi_n(s_0)}{\phi_n(s_1)}}]
      \intersection \product[k < 2][\extensions{\phi_n(s_k) \union
        t_{k,n}}] \subseteq U_{i_n, \differenceset[i_n]{\phi_n(s_0)}
          {\phi_n(s_1)}, n}
  \end{equation*}  
  for all $n \in \N$ and $s_0, s_1 \in \Cantorspace[F_n]$, where
  $\phi_n \from \Cantorspace[F_n] \to \Cantorspace[{\union[m < n]
  [G_m]}]$ is given by $\phi_n(s) = \union[m < n][t_{s(i_m, j_m), m}]$.
  Then the function $\phi \from \Cantorspace[\N \times \N] \to
  \Cantorspace[\N \times \N]$ given by $\support{\phi(c)} \subseteq
  \union[n \in \N][G_n]$ and $\restriction{\phi(c)}{G_n} = t_{c(i_n, j_n),
  n}$ for all $n \in \N$ is continuous.
  
  To see that $\phi$ is a homomorphism from $\equality
  {\Cantorspace}^k \times \Ethree$ to $\equality{\Cantorspace}^k
  \times \Ethree$ for all $k \in \N$, suppose that $c, d \in
  \CantorCantorspace$ are $(\equality{\Cantorspace}^k \times
  \Ethree)$-related, and observe that $\differenceset[k]{t_{c(n), n}}
  {t_{d(n), n}} = \emptyset$ for all $n \in \N$.
  
  To see that $\phi$ is a homomorphism from $\setcomplement
  {\equality{\CantorCantorspace}}$ to $\setcomplement{D}$, note
  that if $c, d \in \CantorCantorspace$ are distinct, then there
  exists $n \in \N$ with the property that $c(i_n, j_n) \neq d(i_n, j_n)$,
  so $\pair{\phi(c)}{\phi(d)} \in U_{i_n, \differenceset[i_n]{\phi_n
  (\restriction{c}{F_n})}{\phi_n(\restriction{d}{F_n})}, n}$, thus $\pair
  {\phi(c)}{\phi(d)} \notin D$.
  
  To see that $\phi$ is a homomorphism from $\setcomplement
  {\Ethree}$ to $\setcomplement{R}$, observe that if $c, d \in
  \CantorCantorspace$ are $\Ethree$-inequivalent, then there
  exists $k \in \N$ such that $\pair{c}{d} \in \D{k} \setminus \D{k+1}$.
  Set $F = \differenceset[k]{c}{d}$, fix $n \in \N$ sufficiently large that
  $F \subseteq F_n$, define $G = \differenceset[k]{\phi_n(\restriction{c}
  {F_n})}{\phi_n(\restriction{d}{F_n})}$, and observe that
  $\differenceset[k]{\phi(c)}{\phi(d)} = G$. As there are arbitrarily large
  $m \ge n$ for which $i_m = k$ and $c(i_m, j_m) \neq d(i_m, j_m)$,
  and therefore $\pair{\phi(c)}{\phi(d)} \in U_{k, G, m}$, it follows that
  $\pair{\phi(c)}{\phi(d)} \notin R$.
\end{propositionproof}

\section{Dichotomies} \label{dichotomies}

We will abuse notation by identifying $\functions{k}{\equality
{\Cantorspace}} \times \Ezero \times \functions{\N}{\equality
{\Cantorspace}}$ with the corresponding equivalence relation on
$\CantorCantorspace$ for all $k \in \N$.

\begin{theorem} \label{dichotomies:general}
  Suppose that $\Gamma$ is a tsi \Polish group, $X$ is an analytic
  metric space, $\Gamma \action X$ is \Borel, $\orbitrelation
  {\Delta}{X}$ is \Borel for all open sets $\Delta \subseteq \Gamma$,
  $\sequence{\Delta_k}[k \in \N]$ is a decreasing sequence of open
  subsets of $\Gamma$ forming a neighborhood basis for $\identity
  {\Gamma}$, and $\Gamma_k$ is the group generated by
  $\Delta_k$. Then exactly one of the following holds:
  \begin{enumerate}
    \item The action $\Gamma \action X$ is $\sigma$-lacunary.
    \item There is a continuous injective homomorphism $\phi \from
      \CantorCantorspace \to X$ from $\sequence{\functions{k}
      {\equality{\Cantorspace}} \times \Ezero \times \functions{\N}
      {\equality{\Cantorspace}}}[k \in \N]$ to $\sequence
      {\orbitequivalencerelation{\Gamma_k}{X}}[k \in \N]$ that is also
      a homomorphism from $\setcomplement{\Ethree}$ to
      $\setcomplement{\orbitequivalencerelation{\Gamma}{X}}$.
  \end{enumerate}
\end{theorem}

\begin{theoremproof}
  Note that condition (2) is equivalent to the apparently weaker
  statement in which $\phi$ is merely \Borel, since we can always
  pass to a dense \Gdelta set $C \subseteq \CantorCantorspace$
  on which $\phi$ is continuous (see, for example, \cite[Theorem
  8.38]{Kechris}), and then compose $\restriction{\phi}{C}$ with the
  map given by Proposition \ref{compositions:continuous}. So by
  \cite[Theorem 5.2.1]{BeckerKechris}, we can assume that
  $\Gamma \action X$ is continuous. 

  By passing to appropriate open subneighborhoods of $\identity
  {\Gamma}$, we can assume that $\Delta_k$ is symmetric and
  $\Delta_{k+1}^2 \subseteq \Delta_k$ for all $k \in \N$. As
  $\Gamma$ is tsi, we can also assume that each $\Delta_k$ is
  conjugation invariant.
  
  Define $G_{i,j} = \orbitrelation{\Delta_i}{X} \setminus \orbitrelation
  {\Delta_j}{X}$ for all $i, j \in \N$. By Propositions \ref
  {lacunary:lacunarytoindependent} and \ref
  {lacunary:independenttolacunary}, condition (1) of Theorems \ref
  {graph:dichotomy} and \ref{dichotomies:general} are equivalent.
  So by Theorem \ref{graph:dichotomy}, it is sufficient to show that
  condition (2) of Theorem \ref{graph:dichotomy} implies condition
  (2) of Theorem \ref{dichotomies:general}. Towards this end,
  suppose that there exist $f \from \N \to \N$ and a continuous
  homomorphism $\phi \from \Cantorspace \to X$ from $\sequence
  {\Gzero[k]}[k \in \N]$ to $\sequence{G_{k, f(k)}}[k \in \N]$.
  
  Appeal to Proposition \ref{compositions:alignment} to obtain a
  continuous homomorphism $\psi \from \Cantorspace \to
  \Cantorspace$ from $\sequence{\Gzero[k]}[k \in \N]$ to $\sequence
  {\Gzero[f^k(0)]}[k \in \N]$. By replacing $\phi$ with $\phi
  \composition \psi$, we can assume that the former is a
  homomorphism from $\sequence{\Gzero[k]}[k \in \N]$ to
  $\sequence{G_{f^k(0), f^{k+1}(0)}}[k \in \N]$. By replacing
  $\sequence{\Delta_k}[k \in \N]$ with $\sequence{\Delta_{f^k(0)}}[k
  \in \N]$, and therefore $\sequence{G_{i,j}}[i,j \in \N]$ with
  $\sequence{G_{f^k(i), f^k(j)}}[i,j \in \N]$, we can assume that $\phi$
  is a homomorphism from $\sequence{\Gzero[k]}[k \in \N]$ to
  $\sequence{G_{k, k+1}}[k \in \N]$.
  
  Fix an enumeration $\sequence{\delta_k}[k \in \N]$ of a countable
  dense subset of $\Gamma$, and for all $k, \ell \in \N$, let $R_{k,
  \ell}$ denote the pullback of $\orbitrelation{\delta_\ell \Delta_k}{X}$
  through $\phi$. Proposition \ref{compositions:focus} then yields
  functions $g_n \from \Cantorspace[<n] \to \N$ for all $n \in \N$ and
  a continuous homomorphism $\psi \from \Cantorspace \to
  \Cantorspace$ from $\sequence{\Gzero[k]}[k \in \N]$ to $\sequence
  {\Gzero[k]}[k \in \N]$ that is also a homomorphism from $\sequence
  {\involutiongraph{s_n}[s]}[n \in \N, s \in \Cantortree]$ to $\sequence
  {R_{k_{n+1+\length{s}}, g_{n+1+\length{s}}(s)}}[n \in \N, s \in
  \Cantortree]$. By replacing $\phi$ with $\phi \composition \psi$ and
  defining $\gamma_{n+1+\length{s}}(s) = \delta_{g_{n+1+\length{s}}
  (s)}$ for all $n \in \N$ and $s \in \Cantortree$, we can assume that
  $\phi$ is also a homomorphism from $\sequence{\involutiongraph
  {s_n}[s]}[n \in \N, s \in \Cantortree]$ to $\sequence{\orbitrelation
  {\gamma_{n+1+\length{s}}(s) \Delta_{k_{n+1+\length{s}}}}{X}}[n \in
  \N, s \in \Cantortree]$.
  
  \begin{lemma} \label{dichotomies:general:composition}
    The function $\phi$ is a homomorphism from $\sequence
    {\involutiongraph{s}}[s \in \Cantortree]$ to $\sequence
    {\orbitequivalencerelation{\Gamma_{k_{\length{s}}}}{X}}[s \in
    \Cantortree]$.
  \end{lemma}
  
  \begin{lemmaproof}
    For each $n \in \N$, let $T_n$ denote the graph on $\Cantorspace
    [n]$ consisting of all pairs of the form $\pair{s_{n - 1 - \length{s}}
    \concatenation \sequence{i} \concatenation s}{s_{n - 1 - \length{s}}
    \concatenation \sequence{1-i} \concatenation s}$, where $i < 2$
    and $s \in \Cantortree[n]$. A simple induction shows that each
    $T_n$ connected.
    
    In particular, it follows that for all $n \in \N$ and $s \in
    \Cantorspace[n]$, there is a $T_n$-path $\sequence{t_\ell}[\ell \le
    m]$ from $s$ to $s_n$. For all $\ell < m$, fix $i_\ell < 2$ and
    $u_\ell \in \Cantortree[n]$ such that $t_\ell = s_{n - 1 - \length
    {u_\ell}} \concatenation \sequence{i_\ell} \concatenation u_\ell$
    and $t_{\ell+1} = s_{n - 1 - \length{u_\ell}} \concatenation
    \sequence{1 - i_\ell} \concatenation u_{\ell}$.
    
    Observe now that if $c \in \Cantorspace$, $i < 2$, and $\ell < m$,
    then $t_\ell \concatenation \sequence{i} \concatenation c$ and
    $t_{\ell+1} \concatenation \sequence{i} \concatenation c$ are
    $\involutiongraph{s_{n - 1 - \length{u_\ell}}}[u_\ell]$-related, so
    $\phi(t_{\ell} \concatenation \sequence{i} \concatenation c)$ and
    $\phi(t_{\ell+1} \concatenation \sequence{i} \concatenation c)$ are
    $\orbitrelation{\gamma_n(u_\ell) \Delta_{k_n}}{X}$-related, thus
    there is an element of $\inverse{(\gamma_n(u_{m-1}) \Delta_{k_n}
    \cdots \gamma_n(u_0) \Delta_{k_n})} \Delta_{k_n}(\gamma_n
    (u_{m-1}) \Delta_{k_n} \cdots \gamma_n(u_0) \Delta_{k_n})$
    sen\-ding $\phi(s \concatenation \sequence{0} \concatenation c)$
    to $\phi(s \concatenation \sequence{1} \concatenation c)$. As the
    conjugation invariance and symmetry of $\Delta_{k_n}$ ensure
    that this product is $\Delta_{k_n}^{2m+1}$, it follows that $\phi(s
    \concatenation \sequence{0} \concatenation c) \mathrel
    {\orbitequivalencerelation{\Gamma_{k_n}}{X}} \phi(s
    \concatenation \sequence{1} \concatenation c)$.
  \end{lemmaproof}
  
  Set $\ell_n = \cardinality{\set{m < n}[k_m = k_n]}$ for all $n \in
  \N$, and define $\psi \from \CantorCantorspace \to \Cantorspace$
  by $\psi(c)(n) = c(k_n, \ell_n)$ for all $c \in \CantorCantorspace$
  and $n \in \N$. Let $D$ and $E$ denote the pullbacks of $\equality
  {X}$ and $\orbitequivalencerelation{\Gamma}{X}$ through $\phi
  \composition \psi$.
  
  \begin{lemma}
    Suppose that $i \in \N$ and $F \in \sets{<\aleph_0}{i \times \N}$.
    Then $E$ is meager in $\D{i}[F]$.
  \end{lemma}
  
  \begin{lemmaproof}
    For all $k \in \N$, let $R_k$ denote the pullback of $\orbitrelation
    {\Delta_k}{X}$ through $\phi \composition \psi$. As $R_{i+2}
    \inverse{R_{i+2}} \subseteq R_{i+1}$, Proposition \ref
    {compositions:meager} ensures that $R_{i+2}$ is meager in
    $\D{i}[F]$. The \Kuratowski-\Ulam theorem therefore ensures
    that for comeagerly-many $c \in \Cantorspace[(i \times \N)
    \setminus F]$ and all $s \in \Cantorspace[F]$, comeagerly-many
    vertical sections of $\set{\pair{d}{d'} \in \Cantorspace[(\N \setminus
    i) \times \N] \times \Cantorspace[(\N \setminus i) \times \N]}[c \union
    s \union d \mathrel{R_{i+2}} c \union \flip{s} \union d']$ are meager,
    so the fact that $\inverse{R_{i+3}} R_{i+3} \subseteq R_{i+2}$
    implies that every vertical section of $\set{\pair{d}{d'} \in
    \Cantorspace[(\N \setminus i) \times \N] \times \Cantorspace[(\N
    \setminus i) \times \N]}[c \union s \union d \mathrel{R_{i+3}} c \union
    \flip{s} \union d']$ is meager. As every vertical section of $\set
    {\pair{d}{d'} \in \Cantorspace[(\N \setminus i) \times \N] \times
    \Cantorspace[(\N \setminus i) \times \N]}[c \union s \union d \mathrel
    {E} c \union \flip{s} \union d']$ is the union of countably-many such
    vertical sections, the \Kuratowski-\Ulam theorem yields that $E$ is
    meager in $\D{i}[F]$.
  \end{lemmaproof}
  
  By composing $\phi \composition \psi$ with the function obtained
  from applying Proposition \ref{compositions:embedding} to $D$ and
  $E$, we obtain the desired homomorphism.
\end{theoremproof}

When $\Gamma$ is non-archimedean, we obtain the following.

\begin{theorem} \label{dichotomies:nonarchimedean}
  Suppose that $\Gamma$ is a non-archimedean tsi \Polish group,
  $X$ is an analytic metric space, $\Gamma \action X$ is \Borel,
  and $\orbitequivalencerelation{\Gamma}{X}$ is \Borel. Then exactly
  one of the following holds:
  \begin{enumerate}
    \item The action $\Gamma \action X$ is $\sigma$-lacunary.
    \item There is a continuous embedding $\pi \from
      \CantorCantorspace \to X$ of $\Ethree$ into
      $\orbitequivalencerelation{\Gamma}{X}$.
  \end{enumerate}
\end{theorem}

\begin{theoremproof}
  By \cite[Theorem 7.1.2]{BeckerKechris}, the orbit equivalence
  relation induced by every open subgroup of $\Gamma$ is \Borel.
  The fact that $\Gamma$ is non-archimedean therefore implies
  that the orbit relation induced by every open subset of $\Gamma$
  is \Borel.

  We can assume that $\Gamma \action X$ is continuous for exactly
  the same reason given at the beginning of the proof of Theorem
  \ref{dichotomies:general}.
  
  Fix a decreasing sequence $\sequence{\Gamma_k}[k \in \N]$ of
  normal subgroups of $\Gamma$ forming a neighborhood basis for
  $\identity{\Gamma}$. In light of Theorem \ref{dichotomies:general},
  we can assume that there is a continuous injective homomorphism
  $\phi \from \CantorCantorspace \to X$ from $\sequence
  {\functions{k}{\equality{\Cantorspace}} \times \Ezero \times
  \functions{\N}{\equality{\Cantorspace}}}[k \in \N]$ to $\sequence
  {\orbitequivalencerelation{\Gamma_k}{X}}[k \in \N]$ that is also
  a homomorphism from $\setcomplement{\Ethree}$ to
  $\setcomplement{\orbitequivalencerelation{\Gamma}{X}}$. But
  the continuity of $\Gamma \action X$ ensures that every such
  function is a reduction of $\Ethree$ to $\orbitequivalencerelation
  {\Gamma}{X}$.
\end{theoremproof}

\begin{acknowledgements}
  I would like to thank Alexander Kechris for pointing out the first
  sentence of the proof of Theorem \ref{dichotomies:nonarchimedean}.
\end{acknowledgements}

\bibliographystyle{amsalpha}
\bibliography{bibliography}

\end{document}